\documentclass[twoside,11pt,leqno]{article}
\usepackage{amsfonts}

 \def\H{{\cal H}}

\def\P{{(\cal {HP})}}

\def\B{B({\cal H})}

\newtheorem{df}{Definition}[section]
\newtheorem{thm}[df]{Theorem} \newtheorem{pro}[df]{Proposition}
 
\newtheorem{rema}[df] {Remark} 
\def\sfstp{{\hskip-1em}{\bf.}{\hskip1em}}

\def\subject#1{\renewcommand{\thefootnote}{}\footnote
{AMS(MOS) subject classification (2010). Primary: {#1}}}

\def\keywords#1{\renewcommand{\thefootnote}{}\footnote
{Keywords: {#1}}}

\def\enddemo{\qed \endtrivlist} \expandafter\let\csname
enddemo*\endcsname=\enddemo

\def\qedsymbol{\ifmmode\bgroup\else$\bgroup\aftergroup$\fi
\vcenter{\hrule\hbox{\vrule
height.5em\kern.5em\vrule}\hrule}\egroup}
\def\qed{\ifmmode\else\unskip\nobreak\fi\quad\qedsymbol}

\title{\bf Expansive operators which are power bounded or algebraic}
\author{\normalsize B.P.~Duggal, I.H.~Kim}
\date{}

\begin{document}

\maketitle \thispagestyle{empty} \vskip-16pt

\subject{ 47A05, 47B47, 47B65; Secondary: 47A55, 47A63, 47A65} \keywords{Expansive/Contractive Hilbert space operator, elementary operator, algebraic operator, power bounded  }
\footnote{The work of the second author was supported by a grant from the National Research Foundation of Korea (NRF), funded by the Korean government (NRF-2019R1F1A1057574).}
\begin{abstract} Given Hilbert space operators $P,T\in B(\H), P\geq 0$ invertible, $T$ is $(m,P)-$ expansive (resp., $(m,P)-$ isometric) for some positive integer $m$ if $\triangle_{T^*,T}^m(P)=\sum_{j=0}^m(-1)^j\left(\begin{array}{clcr}m\\j\end{array}\right){T^*}^jPT^j\leq 0$ (resp., $\triangle_{T^*,T}^m(P)=0$). An  $(m,P)-$ expansive operator $T$ is power bounded if and only if it is a $C_{1\cdot}-$ operator which is similar to an isometry and satisfies $\triangle_{T^*,T}^n(Q)=0$ for some positive invertible operator $Q\in B(\H)$ and all integers $n\geq 1$. If, instead, $T$ is an algebraic $(m,I)-$ expansive operator, then either the spectral radius $r(T)$ of $T$ is greater than one or $T$ is the perturbation of a unitary by a nilpotent such that $T$ is $(2n-1, I)-$ isometric for some positive integers $m_0 \leq m$, $m_0$ odd, and $n \geq \frac{m_0 +1}{2}$.

\end{abstract}


\section {\sfstp Introduction} Let $B(\H)$ denote the algebra of operators, i.e., bounded linear transformations, on an infinite dimensional complex Hilbert space into itself. An operator $T$ is $(m,I)-$ expansive, or simply $m-$ expansive, for some positive integer $m$, if
$$
\triangle_{T^*,T}^m(I)=\sum_{j=0}^m(-1)^j\left(\begin{array}{clcr}m\\j\end{array}\right){T^*}^jT^j\leq 0.
$$
Agler, \cite[Theorem 3,1]{Agl}, characterized subnormality with positivity of $\triangle_{T^*,T}^m(I)$: $\triangle_{T^*,T}^m(I)\geq 0$ if and only if $\|T\| \leq 1$ and $T$ is subnormal. Operators $T$ such that $\triangle_{T^*,T}^m(I)\geq 0$ have been called $m-$ contractive, and operators $T$ such that $\triangle_{T^*,T}^m(I)=0$ are said to be $m-$ isometric \cite{AS}. Classes of $m-$ isometric, $m-$ expansive and $m-$ contractive
 operators have attracted the attention of a large number of authors over the past three or so decades (see \cite{Ath}, \cite{B}, \cite{BJ}, \cite{BMN}, \cite{D}, \cite{EJL}, \cite{Gu},\cite{Gu1}, \cite{JKKL}, \cite{TL} for further references).
For $A,B\in B(\H)$, let $L_A, R_B \in B(B(\H))$ denote respectively the left and the right multiplication operators
$$
L_A(X)=AX \ {\rm and}\ R_B(X)=XB.
$$
Let $\triangle_{A,B} \in B(B(\H))$ denote the elementary operator
$$
\triangle_{A,B}(X))=(I-L_AR_B)(X)=X-AXB.
$$
Then, for positive integers $m$,
$$
\triangle_{A,B}^m(X)=(I-L_AR_B)^m(X)=\sum_{j=0}^m(-1)^j\left(\begin{array}{clcr}m\\j\end{array}\right){A^j}XB^j.
$$
We say in the following that the pair of operators

 $(A,B)\in (m,P)-$ expansive if $\triangle_{A,B}^m(P)\leq 0$;

 $(A,B)\in (m,P)-$ hyperexpansive if $\triangle_{A,B}^t(P)\leq 0$ for all integers $1\leq t\leq m$;

 $(A,B)\in (m,P)-$ contractive if $\triangle_{A,B}^m(P)\geq 0$;

$(A,B)\in (m,P)-$ hypercontractive if $\triangle_{A,B}^t(P)\geq 0$ for all integers $1\leq t\leq m$;

 $(A,B)\in (m,P)-$ isometric if $\triangle_{A,B}^m(P)= 0$.

\

 Recall that an operator $T\in B(\H)$ is power bounded if $\sup_n \|T^n\| \leq M$ for some scalar $M>0$. A well known result says that power bounded $m-$ isometric operators $T$ (i.e., $T$ power bounded and $(T^*,T)\in (m,I)-$ isometric) are isometric; for power bounded pairs $(A,B) \in (m,I)-$ isometric, $A^*$ and $B$ are similar to isometries \cite{DK}. We prove that this result extends to power bounded pairs satisfying the condition that $L_AR_B(Q)\geq 0$ whenever $Q\geq 0$. For power bounded $T\in (m,P)-$ expansive (i.e., $(T^*,T)\in (m,P)-$ expansive) this translates to ``$T$ is a $C_{1\cdot}-$ operator which is similar to an isometry and satisfies $T^*QT=Q$ for some positive invertible operator $Q$". (Thus $T$ is isometric in an equivalent norm: $\|x\|_Q=\langle x,x\rangle_Q^\frac{1}{2}=\|Q^\frac{1}{2}x\|$.)  For operators $T\in (m,P)-$ contractive, it is seen that $T$ is similar to the direct sum of the conjugate of a $C_{0\cdot}-$ contraction with a unitary. Algebraic $(m,P)-$ expansive operators $T$ are not Drazin invertible.  We prove that for such operators $T$ either the spectral radius $r(T)>1$, or, $T$ is the perturbation of a unitary operator by a commuting nilpotent such that $T\in (2n-1)-$ isometric for some integer $n$ (dependent upon $m$). A similar result for algebraic $m-$ contractive operators is not possible.

  \

 We introduce further notation and prove a few complementary results in Section 2, Section 3 in devoted to considering power bounded $(m,P)-$ expansive and $(m,P)-$ contractive operators, and algebraic $m-$ expansive operators are considered in Section 4.

 \section{\sfstp Complementary results}
 Throughout the following $A, B$ and $T$ will denote operators in $B(\H)$, and $P\in B(\H)$ will denote a positive invertible operator. We shall henceforth shorten $(T^*,T)\in (m,P)-\cdots$ to $T\in (m,P)-\cdots$, and $T\in (m,I)-\cdots$ to $T\in m-\cdots$. The spectrum, the approximate point spectrum and the isolated points of the spectrum of $A$ will be denoted by $\sigma(A),\ \sigma_a(A)$ and ${\rm iso}\sigma(A)$, respectively. $T$ is a $C_{0\cdot}-$ operator (resp.,  $C_{1\cdot}-$ operator) if
 $$
 \lim_{n\rightarrow\infty}\|T^nx\|=0 \ {\rm for \ all}\ x \in \H
 $$
 $$
 ({\rm resp.,}\ \inf_{n\in \mathbb{N}}\|T^nx\|>0 \ {\rm for \ all}\ 0 \neq x \in \H );
 $$
$T\in C_{\cdot0}$ if $T^*\in C_{0\cdot}$, $T\in C_{\cdot1}$ if $T^*\in C_{1\cdot}$, and $T\in C_{\alpha\beta}$ if $T\in C_{\alpha\cdot}\cap C_{\cdot\beta} \ (\alpha, \beta=0,1).$ The operator $T$ is weakly $C_{0\cdot}$ (or, weakly stable \cite{Ku}) if $\lim_{n\rightarrow\infty}\langle T^nx,x\rangle =0$ for all $x\in \H$ (equivalently; if $\lim_{n\rightarrow\infty}\langle T^nx,y\rangle =0$ for all $x,y\in \H$). It is well known, \cite{Ker}, that power bounded operators $T$ have an upper triangular representation
$$
T=\left(\begin{array}{clcr}T_1 & T_3\\0&T_2\end{array}\right)\in B(\H_1 \oplus \H_2)
$$
for some decomposition $\H=\H_1\oplus\H_2$ of $\H$ such that $T_1\in C_{0\cdot}$ and $T_2\in C_{1\cdot}.$ Every isometry $V\in B(\H)$ has a direct sum decomposition
$$
V=V_{10}\oplus V_u\in B(\H_c\oplus \H_u), V_{10}\in C_{10} \ {\rm and}\ V_u\in C_{11}
$$
into its completely non-unitary (i.e., unilateral shift) and unitary parts \cite{Ku}.

The following well known result from Douglas \cite{Do} will often be used in the sequel (without further mention).

\begin{thm}\label{thm00} The following statements are pairwise equivalent:
\vskip4pt\noindent (i) $A(\H) \subseteq B(\H)$.

\vskip4pt\noindent (ii) There is a $\mu \geq 0$ such that $AA^* \leq \mu^2BB^*$.

\vskip4pt\noindent (iii) There is an operator $C\in B(\H)$ such that $A=BC$.

If these conditions are satisfied, then the operator $C$ may be chosen so that $\|C\|^2=\inf\{\lambda: AA^*\leq \lambda BB^*\}$, $A^{-1}(0)\subseteq C^{-1}(0)$ and $C(\H)\subseteq B^{-1}(0)^\perp$.
\end{thm}
Suppose that {\em the pair of operators $(A,B)$ preserves order in the sense that $(L_AR_B)(X)\geq 0$ for all $X\in B(\H)$ such that $X\geq 0$}. For all positive integers $n$,
\begin{eqnarray*}
\triangle_{A^n,B^n}^m(P)&=&\left(I-L_{A^n}R_{B^n}\right)^m(P)=(I-L_A^nR_B^n)^m(P)\\
&=& \{L_A^{n-1}\triangle_{A,B}(P)R_B^{n-1}+L_A^{n-2}\triangle_{A,B}(P)R_B^{n-2} +\cdots\\ &+& L_A\triangle_{A,B}(P)R_B+\triangle_{A,B}(P)\}^m\\
&=& \left\{L_A^{n-1}R_B^{n-1}+L_A^{n-2}R_B^{n-2}+\cdots + L_AR_B+I\right\}^m\left(\triangle^m_{A,B}(P)\right).
\end{eqnarray*}
Hence
$$
(A,B)\in (m,P)-{\rm expansive} \Longrightarrow (A^n,B^n)\in (m,P)-{\rm expansive, \ and}
$$
$$
(A,B)\in (m,P)-{\rm contractive} \Longrightarrow (A^n,B^n)\in (m,P)-{\rm contractive}
$$
for all positive integers $n$.

The identity $(a-1)^m=a^m-\sum_{j=0}^{m-1}\left(\begin{array}{clcr}m\\j\end{array}\right)(a-1)^j$ implies
$$
\tilde{\triangle}_{A,B}^m=(L_AR_B-I)^m=(L_AR_B)^m-\sum_{j=0}^{m-1}\left(\begin{array}{clcr}m\\j\end{array}\right)\tilde{\triangle}_{A,B}^{m-1}=(-1)^m\triangle_{A,B}^m.
$$
 If $\tilde{\triangle}_{A,B}^m(P)\leq 0$ for some positive integer $m$, then, since
$$
\tilde{\triangle}_{A,B}^j=L_AR_B\left(\tilde{\triangle}_{A,B}^{j-1}\right)-\tilde{\triangle}_{A,B}^{j-1}
$$
for all integers $j\geq 1$,
\begin{eqnarray*}
\sum_{j=0}^{m-1}\left(\begin{array}{clcr}n\\j\end{array}\right)L_AR_B\left(\tilde{\triangle}_{A,B}^{j}\right)
&=&\sum_{j=0}^{m-1}\left(\begin{array}{clcr}n\\j\end{array}\right)\tilde{\triangle}_{A,B}^{j+1}
+\sum_{j=0}^{m-1}\left(\begin{array}{clcr}n\\j\end{array}\right)\tilde{\triangle}_{A,B}^{j}\\
&=&\left(\begin{array}{clcr}n\\m-1\end{array}\right)\tilde{\triangle}_{A,B}^{m}+\sum_{j=0}^{m-1}\left(\begin{array}{clcr}n+1\\j\end{array}\right)\tilde{\triangle}_{A,B}^{j}.
\end{eqnarray*}
Evidently (see above), $\tilde{\triangle}_{A,B}^{m}(P)\leq 0$ implies
$$
(0\leq)\ (L_AR_B)^m(P)\leq \sum_{j=0}^{m-1}\left(\begin{array}{clcr}m\\j\end{array}\right)\tilde{\triangle}_{A,B}^{j}(P).
$$
We prove
$$
(0\leq)\ (L_AR_B)^n(P)\leq \sum_{j=0}^{m-1}\left(\begin{array}{clcr}n\\j\end{array}\right)\tilde{\triangle}_{A,B}^{j}(P),\ {\rm for \ all}\ n\geq m.
$$
The inequality being true for $n=m$, assume it to be true for $n=t$. Then, since $(A,B)$ preserves order,
\begin{eqnarray*}
(0\leq)\ (L_AR_B)^{t+1}(P)&\leq& \sum_{j=0}^{m-1}\left(\begin{array}{clcr}t\\j\end{array}\right)L_AR_B\left(\tilde{\triangle}_{A,B}^{j}(P)\right)\\
\noindent (1) \hspace{4cm} &=&\left(\begin{array}{clcr}t\\m-1\end{array}\right)\tilde{\triangle}_{A,B}^{m}(P)+\sum_{j=0}^{m-1}\left(\begin{array}{clcr}t+1\\j\end{array}\right)\tilde{\triangle}_{A,B}^{j}(P)\\
&\leq&\sum_{j=0}^{m-1}\left(\begin{array}{clcr}t+1\\j\end{array}\right)\tilde{\triangle}_{A,B}^{j}(P)
\end{eqnarray*}
(since $\tilde{\triangle}_{A,B}^{m}(P)\leq 0$). Thus the inequality is true for $n=t+1$, hence by induction for all integers $n\geq m$.

Observe from (1) that
$$
(0\leq)\ \frac{1}{n^{m-1}}(L_AR_B)^{n}(P)\leq \frac{1}{n^{m-1}}\left\{\left(\begin{array}{clcr}n\\m-1\end{array}\right)\tilde{\triangle}_{A,B}^{m-1}(P)+\sum_{j=0}^{m-2}\left(\begin{array}{clcr}n\\j\end{array}\right)
\tilde{\triangle}_{A,B}^{j}(P)\right\}.
$$
Since $\left(\begin{array}{clcr}n\\m-1\end{array}\right)$ is of the order of $n^{m-1}$ and $\left(\begin{array}{clcr}n\\m-2\end{array}\right)$ is of the order of $n^{m-2}$ for large $n$, letting $n\rightarrow \infty$ we have
$$
0\leq \tilde{\triangle}_{A,B}^{m-1}(P) \ \left( \Longleftrightarrow (-1)^m\triangle_{A,B}^{m-1}(P)\leq 0\right).
$$
In conclusion, we have:

\begin{pro}\label{pro00} If the pair $(A,B)$ preserves order, then
\vskip4pt\noindent (i) $m$ positive even and $(A,B)\in (m,P)-$ expansive implies $ (A,B)\in (m-1,P)-$ expansive ;

\vskip4pt\noindent (ii) $m$ positive odd and $(A,B)\in (m,P)-$ contractive implies $(A,B)\in (m-1,P)-$ contractive.
\end{pro}

For pairs $(T^*,T)$ this translates to ({\em cf} \cite{Gu1}):

\begin{pro}\label{pro01} If $T\in (m,P)-$ expansive for some even positive integer $m$ (resp., $T\in (m,P)-$ contractive for some odd positive integer $m$), then $T\in (m-1,P)-$ expansive (resp., $T\in (m-1,P)-$ contractive).
\end{pro}

\

\section{\sfstp Power bounded operators}
Proposition \ref{pro00} does not extend to odd positive integers $m$  for $(m,P)-$ expansive (resp, even positive integers $m$ for $(m,P)-$ contractive) operators $T$: for if it were so, then one would have that $T\in (m,P)-$ expansive implies $T\in (m,P)-$ hyperexpansive (resp., $T\in (m,P)-$ contractive implies $T\in (m,P)-$ hypercontractive). A class of operators where Proposition \ref{pro00} does have an extension to all $m$ is that of power bounded operators. We have:

\begin{thm}\label{thm30} If $A,B$ are power bounded, the pair $(A,B)$ preserves order \newline and $(A,B)\in (m,P)-$ expansive (resp., $(A,B)\in (m,P)-$ contractive), then    $(A,B)\in (m,P)-$ hyperexpansive  (resp., $(A,B)\in (m,P)-$ hypercontractive).
\end{thm}

\begin{demo}
In view of Proposition \ref{pro00}, we have only to prove that $m$ odd, $(A,B)\in (m,P)-$ expansive implies $(A,B)\in (m-1,P)-$ expansive and $m$ even, $(A,B)\in (m,P)-$ contractive implies $(A,B)\in (m-1,P)-$ contractive. And for this it is sufficient to prove that
$$
\tilde{\triangle}_{A,B}^{m}(P)\geq 0\Longrightarrow \tilde{\triangle}_{A,B}^{m-1}(P)\leq 0,
$$
since by definition
$$
\triangle_{A,B}^{m}(P)\leq 0\Longleftrightarrow \tilde{\triangle}_{A,B}^{m}(P)\geq 0,\ m \ {\rm odd}
$$
and
$$
\triangle_{A,B}^{m}(P)\geq 0 \Longleftrightarrow \tilde{\triangle}_{A,B}^{m}(P)\geq 0,\ m \ {\rm even}.
$$
If $\tilde{\triangle}_{A,B}^{m}(P)\geq 0$, then
$$
\tilde{\triangle}_{A,B}^{m}(P)=(L_AR_B)^m(P)-\sum_{j=0}^{m-1}\left(\begin{array}{clcr}m\\j\end{array}\right)\tilde{\triangle}_{A,B}^{j}(P)\geq 0.
$$
By hypothesis, $(A,B)$ preserves order. Hence, since
\begin{eqnarray*}
&&(L_AR_B)\left\{(L_AR_B)^t-\sum_{j=0}^{m-1}\left(\begin{array}{clcr}t\\j\end{array}\right)\tilde{\triangle}_{A,B}^{j}\right\}\\
&=&(L_AR_B)^{t+1}-\left\{\sum_{j=0}^{m-1}\left(\begin{array}{clcr}t+1\\j\end{array}\right)\tilde{\triangle}_{A,B}^{j}
+\left(\begin{array}{clcr}t\\m-1\end{array}\right)\tilde{\triangle}_{A,B}^{m}\right\},\\
\end{eqnarray*}
an induction argument shows that
\begin{eqnarray*}
\noindent (2) \hspace{2cm}0&\leq& (L_AR_B)^{n}(P)-\left\{\sum_{j=0}^{m-1}\left(\begin{array}{clcr}n\\j\end{array}\right)\tilde{\triangle}_{A,B}^{j}(P)
+\left(\begin{array}{clcr}n-1\\m-1\end{array}\right)\tilde{\triangle}_{A,B}^{m}(P)\right\}\\
&\leq& (L_AR_B)^{n}(P)-\sum_{j=0}^{m-1}\left(\begin{array}{clcr}n\\j\end{array}\right)\tilde{\triangle}_{A,B}^{j}(P)
\end{eqnarray*}
for all integer $n\geq m$.

The power bounded hypothesis on $A,B$ implies
$$
\left| \langle (L_AR_B)^{n}(P)x,x\rangle \right| \leq \left\|P^\frac{1}{2}\right\|^2\left\|A^{*n}\right\|\left\|B^n\right\|\|x\|^2\leq M\|x\|^2
$$
for some scalar $M>0$. Hence, since
$$
\sum_{j=0}^{m-1}\left(\begin{array}{clcr}n\\j\end{array}\right)\tilde{\triangle}_{A,B}^{j}=\left(\begin{array}{clcr}n\\m-1\end{array}\right)\tilde{\triangle}_{A,B}^{m-1}
+\sum_{j=0}^{m-2}\left(\begin{array}{clcr}n\\j\end{array}\right)\tilde{\triangle}_{A,B}^{j},
$$
$\left(\begin{array}{clcr}n\\m-1\end{array}\right)$ is of the order of $n^{m-1}$ and $\left(\begin{array}{clcr}n\\j\end{array}\right)$ is of the order of $n^{m-2}$ $(0\leq j\leq m-2)$ as $n\rightarrow \infty$, it follows upon dividing the inequality in (2) by $n^{m-1}$ and letting $n\rightarrow \infty$ that
$$
-\tilde{\triangle}_{A,B}^{m-1}(P)\geq 0\Longleftrightarrow \tilde{\triangle}_{A,B}^{m-1}(P)\leq 0.
$$
\end{demo}

The following theorem says that for power bounded order preserving pairs of operators $(A,B)\in (m,P)-$expansive, $A$ and $B$ have a simple form: $B$ is similar to an isometry and $A$ is similar to a co-isometry.

\begin{thm}\label{thm31} Given power bounded operators $A,B$ such that $(A,B)$ preserves order, if $(A,B)\in (m,P)-$ expansive, then there exist positive operators $P_i$ and isometries $V_i,\ i=1,2$, such that $A=P_1^{-1}V_1^*P_1$ and $B=P_2^{-1}V_2P_2$.
\end{thm}

\begin{demo} Since $\triangle_{A,B}^{m}(P)\leq 0$ implies $\triangle_{A^n,B^n}^{m}(P)\leq 0$ for all positive integers $n$, we have:
\begin{eqnarray*}
&&(A,B)\in (m,P)-{\rm expansive} \Longrightarrow \triangle_{A^n,B^n}^{m}(P)\leq 0\\
&\Longleftrightarrow& P\leq \sum_{j=1}^{m}(-1)^{j+1}\left(\begin{array}{clcr}m\\j\end{array}\right)A^{nj}PB^{nj}\\
&\Longleftrightarrow& I\leq \sum_{j=1}^{m}(-1)^{j+1}\left(\begin{array}{clcr}m\\j\end{array}\right)
\left(P^{-\frac{1}{2}}A^nP^{\frac{1}{2}}\right)^j\left(P^{\frac{1}{2}}B^nP^{-\frac{1}{2}}\right)^{j-1}\left(P^{\frac{1}{2}}B^nP^{-\frac{1}{2}}\right)\\
&\Longrightarrow& \|x\|\leq \left\|\sum_{j=1}^{m}(-1)^{j+1}\left(\begin{array}{clcr}m\\j\end{array}\right)
\left(P^{-\frac{1}{2}}A^n(P^{\frac{1}{2}}\right)^j\left(P^{\frac{1}{2}}B^n(P^{-\frac{1}{2}}\right)^{j-1}\right\|\left\|P^{\frac{1}{2}}B^nP^{-\frac{1}{2}}x\right\|\\
&\Longrightarrow& \|x\|\leq M_0 \left\|P^{\frac{1}{2}}B^nP^{-\frac{1}{2}}x\right\|
\end{eqnarray*}
for some scalar $M_0>0$ and all $x\in \H$. The operator $P^{\frac{1}{2}}BP^{-\frac{1}{2}}$ being power bounded, there exists a scalar $M_1>0$ such that
$$
\frac{1}{M_0} \|x\|\leq \left\|\left(P^{\frac{1}{2}}BP^{-\frac{1}{2}}\right)^nx\right\| \leq M_1\|x\|
$$
for all $x\in \H$. Hence there exists an invertible operator $S$ and an isometry $V$ such that
$$
P^{\frac{1}{2}}BP^{-\frac{1}{2}}=S^{-1}VS \Longleftrightarrow B=\left(SP^{\frac{1}{2}}\right)^{-1}V\left(SP^{\frac{1}{2}}\right)
$$
\cite{KR}. But then
$$
B^*P^{\frac{1}{2}}S^*SP^{\frac{1}{2}}B=P^{\frac{1}{2}}S^*SP^{\frac{1}{2}} \Longleftrightarrow B^*P_1^2B=P_1^2
$$
for some invertible positive operator $P_1^2=P^{\frac{1}{2}}S^*SP^{\frac{1}{2}}$.

\noindent Conclusion: there exists an isometry $V_1$ and a positive invertible operator $P_1$ such that
$$
B^*P_1=P_1V_1^*\Longleftrightarrow B=P_1^{-1}V_1P_1.
$$
To complete the proof, we apply the above argument to
$$
\triangle_{B^*,A^*}^{m}(P)\leq 0 \left(\Longleftrightarrow \triangle_{A,B}^{m}(P)\leq 0 \right)
$$
to conclude the existence of an invertible positive operator $P_2$ and an isometry $V_2$ such that $A=P_2^{-1}V_2^*P_2$.
\end{demo}

For $(m,P)-$ contractive pairs $(A,B)$ of power bounded operators Theorem \ref{thm30} implies
$$
\triangle_{A,B}(P)\geq 0 \Longleftrightarrow \left(P^{-\frac{1}{2}}A(P^{\frac{1}{2}}\right)\left(P^{\frac{1}{2}}B(P^{-\frac{1}{2}}\right)\leq I.
$$
Letting $A=B^*=T^*$, it then follows that:

if  $T\in (m,P)-$ expansive (resp., $T\in (m,P)-$ contractive), then $T$ is similar to an isometry (resp., $T$ is similar to a contraction, hence similar to a part of a co-isometric operator \cite[Lemma 7.1]{Ku}).

More is true.

\begin{thm}\label{thm32}  The following conditions are pairwise equivalent for $(m,P)$-expansive operators  $T\in B(\H)$.

\vskip4pt\noindent (i) $T$ is power bounded.

\vskip4pt\noindent (ii) $T$ is (a $C_{1\cdot}$-operator which is) similar to an isometry.

\vskip4pt\noindent (iii) There exists a positive invertible operator $Q$ such that $T\in (n,Q)$-isometric for all integers $n\geq 1$.

\vskip4pt\noindent (iv) There exists a positive invertible operator $Q$  and an equivalent norm $||.||_Q$ on $\H$ induced by the inner product $\langle . , . \rangle_Q=\langle Q., .\rangle$  such that $T$ is $n$-isometric for all integers $n\geq 1$ in this new norm.
\end{thm}

\begin{demo} $(i)\Longrightarrow (ii).$  If $T\in (m,P)-$ expansive, then (see above) there exists a positive invertible operator $P_1\in B(\H)$ and an isometry $V_1\in B(\H)$ such that $P_1T=V_1P_1$. The operator $T$ being power bounded, there exists a direct sum decomposition $\H=\H_{11}\oplus \H_{12}$ of $\H$ such that
$$
T=\left(\begin{array}{clcr}T_1 & T_3\\0&T_2\end{array}\right)\in B(\H_{11} \oplus \H_{12}), \ T_1\in C_{0\cdot} \ {\rm and} \ T_2\in C_{1\cdot}
$$
\cite{Ker}. Decompose $V_1$ into its completely non-unitary (i.e., forward unilateral shift) and unitary parts by
$$
V_1=V_{10}\oplus V_{1u} \in B(\H_{10} \oplus \H_{20}).
$$
Let $P_1\in  B(\H_{11} \oplus \H_{12}, \H_{10} \oplus \H_{20})$ have the matrix representation
$$
P_1=\left(\begin{array}{clcr}P_{11} & P_{12}\\P_{12}^*&P_{22}\end{array}\right).
$$
Then
\begin{eqnarray*}
P_1T&=&V_1P_1 \Longrightarrow P_{12}^*T_1=V_{1u}P_{12}^* \Longrightarrow P_{12}^*T_1^n=V_{1u}^nP_{12}^* \ ({\rm all \ positive\ integers}\ n)\\
&\Longrightarrow& \left\|P_{12}^*x\right\|=\left\|V_{1u}^nP_{12}^*x\right\|=\left\|P_{12}^*T_1^nx\right\|\leq \left\|P_{12}^*\right\|\left\|T_1^nx\right\|
\end{eqnarray*}
for all $x\in \H_{11}$. Since $T_1\in C_{0\cdot}$,
$$
\left\|P_{12}^*x\right\|\rightarrow 0 \ {\rm as}\ n \rightarrow \infty \Longleftrightarrow P_{12}^*=0.
$$
Hence
$$
P=P_{11}\oplus P_{22}, \ P_{11}\ {\rm and}\ P_{22}\geq 0 \ {\rm invertible},
$$
and
$$
P_1T=V_1P_1 \Longrightarrow P_{11}T_3=0,\ P_{11}A_1=V_{10}P_{11}.
$$
Consequently, $T_3=0$ and
\begin{eqnarray*}
&&P_{11}A_1=V_{10}P_{11} \Longrightarrow P_{11}A_1^n=V_{10}^nP_{11} \ ({\rm all \ positive\ integers}\ n)\\
&\Longrightarrow& \left\|P_{11}x\right\|=\left\|V_{10}^nP_{11}x\right\|=\left\|P_{11}A_1^nx\right\|\leq \left\|P_{11}\right\|\left\|A_1^nx\right\|\rightarrow 0 \ {\rm as}\ n \rightarrow \infty\\
&&\ \ \ \ \ \ \ \ ({\rm since }\ A_1\in C_{0\cdot})\\
&\Longrightarrow& \left\|P_{11}x\right\|=0 \Longleftrightarrow P_{11}=0 \ {\rm or}\ x=0.
\end{eqnarray*}
Since $P_{11}$ is invertible, we must have $\H_{11}=\{0\}$, and then $T$ is a $C_{1\cdot}$-operator such that $T=P_1^{-1}V_1P_1$ .

\vskip4pt $(ii)\Longrightarrow (iii).$  Evident, since $(ii)$ holds implies 
$$ T=P^{-1}_1VP_1 \Longrightarrow T^*QT=Q, Q=P^{2}_1\Longrightarrow T\in (n,Q)-{\rm isometric}$$
for all positive integers $n\geq 1$.

\vskip4pt $(iii)\Longrightarrow (iv).$  The operator $Q\geq 0$ being invetible, $||.||_Q$ is an equivalent norm on $\H$ \cite{KD} such that $\sum_{j=0}^n{(-1)^j\left(\begin{array}{clcr}n\\j\end{array}\right)||T^jx||^2_Q}=0$ for  integers $n\geq 1$  and all $x\in\H$.

\vskip4pt $(iv)\Longrightarrow (i).$  Evident, since $T\in (n,Q)$-isometric implies $T^p\in (n,Q)$-isometric for all integers $p\geq 1$, in particular
\begin{eqnarray*} & & 0=||x||^2_Q-||T^px||^2_Q=\langle (Q-T^{*p}QT^p)x,x \rangle  \  {\rm for \  all} \   x\in\H  \Longleftrightarrow Q=T^{*p}QT^p\\&\Longleftrightarrow&
 {\rm there \  exists \ an \ isometry} \ V \ {\rm such \ that} \ T^{*p}Q^{\frac{1}{2}}=Q^{\frac{1}{2}}V^*\Longleftrightarrow T^p=Q^{-\frac{1}{2}}VQ^{\frac{1}{2}}\\
&\Longrightarrow& \sup_p||T^p||\leq ||Q^{-\frac{1}{2}}||||Q^{\frac{1}{2}}||<\infty.\end{eqnarray*}
This completes the proof. 
\end{demo}
For $(m,P)$-contractive power bounded operators, we have:

\begin{thm}\label{thm310} If $T$ is a power bounded $(m,P)$-contractive operator in $\B$, then $T$ is similar to the direct sum of the adjoint of a $C_{0\cdot}$-contraction with a unitary.
\end{thm}
\begin{demo} If $T\in (m,P)-$contractive is power bounded, then $T\in (m,P)-$hypercontractive (by Theorem \ref{thm00}) and hence
$$
\triangle_{T^*,T}(P)\geq 0 \Longleftrightarrow P\geq T^*PT.
$$
Consequently, there exists a contraction $C\in B(\H)$ such that
$$
P^\frac{1}{2}C=T^*P^\frac{1}{2}.
$$
The contraction $C$ has a decomposition, the Foguel decomposition \cite{Ku},
$$
C=Z\oplus U\in B(\H_c\oplus \H_c^\perp),
$$
$$
\H_e=\{x\in\H: \langle C^nx,y\rangle \rightarrow 0 \ {\rm as} \ n\rightarrow \infty,\ {\rm all}\ y\in\H\},
$$
where $U$ is unitary and
$$
\lim_{n\rightarrow \infty}\langle Z^nx,x\rangle=0 \ {\rm for\ all}\ x\in \H_c
$$
(i.e., $Z\in B(\H_c)$ is weakly $C_{0\cdot}$). Letting, as before
$$
T=\left(\begin{array}{clcr}T_1 & T_3\\0&T_2\end{array}\right)\in B(\H_{11} \oplus \H_{12}), \ T_1\in C_{0\cdot} \ {\rm and} \ T_2\in C_{1\cdot},
$$
and letting
$$
P^\frac{1}{2}=\left(\begin{array}{clcr}P_{11} & P_{12}\\P_{12}^*&P_{22}\end{array}\right)\in B(\H_{11} \oplus \H_{12}, \H_c\oplus \H_c^\perp).
$$
the equality
\begin{eqnarray*}
&&P_{12}U=T_1^*P_{12} \Longleftrightarrow U^*P_{12}^*=P_{12}^*T_1 \Longrightarrow U^{*n}P_{12}^*=P_{12}^*T_1^n\\
&\Longrightarrow& \left\|P_{12}^*x\right\|=\left\|U^{*n}P_{12}^*x\right\|=\left\|P_{12}^*T_1^nx\right\|\leq \left\|P_{12}^*\right\|\left\|T_1^nx\right\|\ ({\rm all} \ x\in \H_{11})\\
&\Longrightarrow& \left\|P_{12}^*x\right\|\leq \left\|P_{12}^*\right\|\lim_{n\rightarrow \infty}\left\|T_1^nx\right\|=0\\
&\Longrightarrow& P_{12}=0, \ P^\frac{1}{2}=P_{11}\oplus P_{22},\ P_{11} \ {\rm and}\ P_{22}\geq 0 \ {\rm invertible}.
\end{eqnarray*}
Considering now  $T_3^*P_{11}=0$ it follows that
$$
T_3=0, \ T=T_1\oplus T_2, \ T_1^*=P_{11}ZP^{-1}_{11},\ T^*_2=P_{22}UP_{22}^{-1}
$$
and $T$ is  similar to the direct sum of the adjoint of a  $C_0\cdot-$ contraction (hence, a weakly $C_{0\cdot}-$ contraction)  with a unitary.
\end{demo}

\

It is clear from the above that in the case in which $T\in (m,P)-$ isometric, then ($T\in (m,P)-$ expansive $\wedge$ $(m,P)-$ contractive) $P_1-T^*P_1T=0$, where the similarity $P_1$ may be chosen to be the operator $P$. In particular, if $P=I$, then $T$ is isometric.

\

\section{\sfstp Algebraic $T$}
If $T\in B(\H)$ is an algebraic operator (i.e., there exists a polynomial $q$ such that $q(T)=0$), then $T$ has a representation
$$
T=\bigoplus_{i=1}^tT\mid_{\H_0(T-\lambda_i I)},\ \H=\bigoplus_{i=1}^t{\H_0(T-\lambda_i I)}
$$
for some positive integer $t$ and scalars $\lambda_i$, where
\begin{eqnarray*}
\H_0(T-\lambda_iI)&=&\left\{x\in \H : \lim_{n\rightarrow \infty}\left\|(T-\lambda_iI)^nx\right\|^\frac{1}{n}=0\right\}\\
&=& (T-\lambda_iI)^{-p_i}(0)
\end{eqnarray*}
for some positive integer $p_i$. The points $\lambda_i$ are poles of the resolvent of $T$ of order $p_i$ and (therefore) each $T_i=T\mid_{\H_0(T-\lambda_i I)}$ has a representation
$$
T_i=\lambda_iI_i+N_i, \ 1\leq i\leq t,
$$
where $I_i$ is the identity of $B(\H_0(T-\lambda_iI))$ and $N_i$ is $p_i$-nilpotent. Evidently,
$$
T=\bigoplus_{i=1}^tT_i = \bigoplus_{i=1}^t(\lambda_i I_i+N_i)=T_0+N,
$$
where $T_0$ is a normal operator with
$$
\sigma(T_0)=\sigma_a(T_0)=\sigma(T)=\{\lambda_i\}_{i=1}^t
$$
and $N$ is a nilpotent of order $p=\max\{p_i: 1\leq i\leq t\}$.

Assume now that $T\in (m,P)-$ expansive, $P\geq 0$ invertible (as before). If $\lambda \in \sigma_a(T)$, then there exists a sequence of unit vectors $\{x_n\}\subset\H$ such that $\lim_{n\rightarrow \infty}\left\|(T-\lambda I)x_n\right\|=0$ and
\begin{eqnarray*}
\lim_{n\rightarrow \infty}\langle \triangle_{T^*,T}^m(P)x_n,x_n\rangle&=&\lim_{n\rightarrow \infty}\sum_{j=0}^{m}(-1)^{j}\left(\begin{array}{clcr}m\\j\end{array}\right)\left\|P^\frac{1}{2}T^jx\right\|^2\\
&=&\lim_{n\rightarrow \infty}\sum_{j=0}^{m}(-1)^{j}\left(\begin{array}{clcr}m\\j\end{array}\right)| \lambda |^{2j}\left\|P^\frac{1}{2}x_n\right\|^2\\
&=&\lim_{n\rightarrow \infty}\left(1-| \lambda |^{2}\right)^m\left\|P^\frac{1}{2}x_n\right\|^2\leq 0.
\end{eqnarray*}
Since $P\geq 0$ is invertible, we must have
$$
|\lambda|=1 \ {\rm if}\ m\ {\rm is \ even} \ (\Longrightarrow \sigma_a(T) \subseteq \partial\mathbb{D} \ {\rm if}\ m\ {\rm is \ even})
$$
and
$$
|\lambda|\geq 1 \ {\rm if}\ m\ {\rm is \ odd} \ (\Longrightarrow \sigma_a(T) \subseteq \mathbb{C} \setminus\mathbb{D} \ {\rm if}\ m\ {\rm is \ odd})
$$

Algebraic $(m,P)-$ expansive operators cannot be Drazin invertible (hence are invertible). To see this, let $T$ be an algebraic $(m,P)-$ expansive Drazin invertible operator. Then there exists a decomposition $\H=\H_1\oplus\H_2$ of $\H$, a decomposition $T=T\mid_{\H_1}\oplus T\mid_{\H_2}=T_1\oplus T_2$ of $T$ such that $T_1$ is invertible and $T_2$ is $p-$ nilpotent for some positive integer $p$. Since
$$
T\in (m,P)-{\rm expansive}\Longrightarrow T^p\in (m,P)-{\rm expansive},
$$
letting $P\in B(\H_1\oplus\H_2)$ have the representation $P=\left[P_{ik}\right]_{i,k=1}^2$, we have
\begin{eqnarray*}
0&\geq& \sum_{j=0}^{m}(-1)^{j}\left(\begin{array}{clcr}m\\j\end{array}\right)T^{*p_j}PT^{p_j}\\
&=&\left[\sum_{j=0}^{m}(-1)^{j}\left(\begin{array}{clcr}m\\j\end{array}\right)T_i^{*p_j}P_{ik}T_k^{p_j}\right]_{i,k=1}^2\\
&=&\left(\begin{array}{clcr}\sum_{j=0}^{m}(-1)^{j}\left(\begin{array}{clcr}m\\j\end{array}\right)T_1^{*p_j}P_{11}T_1^{p_j}&P_{12}\\P_{21}&P_{22}\end{array}\right)\\
&\Longrightarrow& P_{22}=0 \ ({\rm since}\ P_{22}\geq 0).
\end{eqnarray*}
But then, by the positivity of $P, \  P_{12}=P_{21}^*=0$ (\cite{A}, Theorem I.1). Since $P$ is invertible, this is a contradiction.

The following theorem says that for an algebraic $(m,P)-$ expansive operator $T$, either $r(T)>1$ or $T$ is the direct sum of a unitary with a nilpotent $(2n-1)-$ isometric operator for some positive integer $n$.

\begin{thm}\label{thm10} If $T\in B(\H)$ is an algebraic $m-$ expansive operator such that $r(T) \leq 1$, then:
\vskip4pt\noindent (i) $T$ is a perturbation of a unitary by a commuting nilpotent;

\vskip4pt\noindent (ii) there exist positive integers $m_0$ and $n$,
$$
m_0\leq m,\ m_0 \ {\rm odd}, \ n\geq \frac{m_0+1}{2},
$$
such that
$$
T\in (2n-1)-{\rm isometric}.
$$
\end{thm}

\begin{demo} We consider $m$ even and $m$ odd cases separately. If $m$ is even, then (as seen above)
$$
\sigma(T)=\sigma_a(T)\subseteq \partial\mathbb{D} \Longrightarrow \sigma(T_0) \subseteq \partial\mathbb{D}
$$
hence the normal operator $T_0$ is a unitary (and $T=T_0+N$, $ [T_0,N]=0$, is the perturbation of $T_0$ by a nilpotent).  Then
\begin{eqnarray*}
\triangle_{T^*,T}&=&\left(I-L_{T^*}R_T\right)=\left(I-L_{T_0^*+N^*}R_{T_0+N}\right)\\
&=&\left(I-L_{T_0^*}R_{T_0}\right)-\left\{L_{N^*}R_{T_0}+L_{T_0^*+N^*}R_N\right\}\\
&=&\triangle_{T_0^*,T_0}-\left\{L_{N^*}R_{T_0}+L_{T_0^*+N^*}R_{N}\right\}\\
&=&\triangle_{T_0^*,T_0}-S \ ({\rm say})
\end{eqnarray*}
and
\begin{eqnarray*}
\triangle_{T^*,T}^t(I)&=&\left(\sum_{j=0}^{t}(-1)^j\left(\begin{array}{clcr}t\\j\end{array}\right)\triangle_{T_0^*,T_0}^{t-j}S^{t}\right)(I)\\
&=&\sum_{j=0}^{t}(-1)^{j}\left(\begin{array}{clcr}t\\j\end{array}\right)S^{j}\triangle_{T_0^*,T_0}^{t-j}(I)
\end{eqnarray*}
(since $[T_0,N]=0$). Evidently,
$$
T_0\in 1-{\rm isometric}  \  (\Longleftrightarrow \triangle_{T_0^*,T_0}(I)=0);
$$
hence
$$
\triangle_{T^*,T}^t(I)= (-1)^tS^t=(-1)^t \left(\sum_{k=0}^{t}\left(\begin{array}{clcr}t\\k\end{array}\right)R_{T^*_0}^{t-k}L^k_{T_0^*+N^*}L_{N^*}^{t-k}R_N^k\right)(I)
$$
This implies that if $N$ is $n-$ nilpotent and $t=2n-1$, then $S=0$ and, consequently, $T\in (2n-1)$-isometric. We prove that $n\geq \frac{m_0+1}{2}$. By hypothesis (above) the odd integer $m_0$ is the smallest positive integer such that
\begin{eqnarray*}
&&\langle \triangle_{T^*,T}^{m_0-1}(I)x_0,x_0\rangle=\sum_{j=0}^{m_0-1}(-1)^j\left(\begin{array}{clcr}m_0-1\\j\end{array}\right)\left\|T^jx_0\right\|^2>0\\
&\Longleftrightarrow& \langle S^{m_0-1}x_0,x_0\rangle=\sum_{k=0}^{m_0-1}(-1)^j\left(\begin{array}{clcr}m_0-1\\k\end{array}\right){N^*}^{m_0-1-k}\langle {T_0^*+N^*}^{m_0-1}T_0^{m_0-1-1}x_0,x_0\rangle N^k >0.
\end{eqnarray*}
Since $n=\frac{m_0-1}{2}$ forces
$$
\langle S^{m_0-1}x_0,x_0\rangle=0,
$$
we must have $N^n\neq 0$ for all $n\leq \frac{m_0-1}{2}$.

\

\noindent If $m$ is odd, then
$$
\sigma(T)=\sigma_a(T)\subseteq \mathbb{C}\setminus\mathbb{D}
$$
and the spectral radius
$$
r(T)=\max\{|\lambda|: \lambda\in \sigma(T)\}
$$
satisfies $r(T)=1$ or $r(T)>1$. If  $r(T)=1$, then $\sigma(T)=\sigma_a(T)\subseteq \partial\mathbb{D} $ and $T=T_0+N$ is the perturbation of a unitary by a commuting nilpotent. The argument above applies, and the proof follows.
\end{demo}

The theorem fails in the case in which $m$ is odd and $r(T)>1$. Consider, for example, the operator $T=\alpha I$, where $|\alpha|>1$. Then
$$
\triangle_{T^*,T}^{2m+1}(I)=\sum_{j=0}^{2m+1}(-1)^j\left(\begin{array}{clcr}2m+1\\j\end{array}\right)|\alpha|^{2j}=\left(1-|\alpha|^{2}\right)^{2m+1}<0
$$
for all integers $m\geq 0$. Observe here that
$$
\tilde{\triangle}_{T^*,T}^m(I)>0
$$
for all positive integers $m$ (i.e., the operator $T$ is $m$-alternatingly expansive \cite[Definition 1.1(7)]{Gu}). Is this typical of operators $T\in m-$ expansive for some odd positive integer $m$ with $r(T)>1$ ? The operator $T$ of the example evidently satisfies $T^*T>I$: The following proposition proves that invertible operators $T$ such that $T\in (m,P)-$ expansive, $P\geq 0$ invertible and $T^*T\geq 1$ are indeed $(m,P)-$ alternatingly expansive.

\begin{pro}\label{pro11} If an invertible operator $T\in (m,P)-$ expansive, $P\geq 0$ invertible, satisfies $T^*T\geq 1$, then $T\in (m,P)-$ alternatingly expansive.
\end{pro}

\begin{demo} The hypotheses imply that $T^{-1}$ is a contraction, hence power bounded, such that
$$
\tilde{\triangle}_{{T^*}^{-1},T^{-1}}^m(P)\leq 0.
$$
Consequently,
$$
\tilde{\triangle}_{{T^*}^{-1},T^{-1}}^m(P)=\left\{ \begin{array}{l}
\triangle_{{T^*}^{-1},T^{-1}}^m(P)\leq 0 \ {\rm if} \ m \ {\rm is\ even}\\
\triangle_{{T^*}^{-1},T^{-1}}^m(P)\geq 0 \ {\rm if} \ m \ {\rm is\ odd},
\end{array}\right.
$$
and this (by Proposition \ref{pro01}) implies
$$
T^{-1}\in \left\{ \begin{array}{l}
(m,P)-{\rm hyperexpansive \ if} \ m \ {\rm is\ even}\\
(m,P)-{\rm hypercontractive \ if} \ m \ {\rm is\ odd}.
\end{array}\right.
$$
Since
$$
\triangle_{{T^*}^{-1},T^{-1}}^t(P)\leq 0\Longrightarrow \left\{ \begin{array}{l}
\triangle_{{T^*},T}^t(P)\leq 0 \ {\rm if} \ t \ {\rm is\ even}\\
\triangle_{{T^*},T}^t(P)\geq 0 \ {\rm if} \ t \ {\rm is\ odd}
\end{array}\right.
$$
and
$$
\triangle_{{T^*}^{-1},T^{-1}}^t(P)\geq 0 \Longrightarrow \left\{ \begin{array}{l}
\triangle_{{T^*},T}^t(P)\geq 0 \ {\rm if} \ t \ {\rm is\ even}\\
\triangle_{{T^*},T}^t(P)\leq 0 \ {\rm if} \ t \ {\rm is\ odd},
\end{array}\right.
$$
the proof follows
\end{demo}

\begin{rema} We remark in closing that a similar analysis does not hold for $(m,P)-$ contractive algebraic operators. Thus $T=\alpha I \oplus 0 \in B(\H\oplus\H)$ is Drazin invertible $(m, P_1\oplus P_2)-$ contractive operator, $P_1$ and $P_2\in B(\H)$ are positive invertible, for all scalars $\alpha$ if $m$ is even and for scalars $\alpha$ such that $|\alpha|\leq 1$ if $m$ is odd.
\end{rema}


\vskip10pt \noindent\normalsize\rm B.P. Duggal, Faculty of Mathematics, Visgradska 33, 1800 Ni\v{s}, Serbia.\\
\noindent\normalsize \tt e-mail: bpduggal@yahoo.co.uk

\vskip6pt\noindent \noindent\normalsize\rm I. H. Kim, Department of
Mathematics, Incheon National University, Incheon, 22012, Korea.\\
\noindent\normalsize \tt e-mail: ihkim@inu.ac.kr

\end{document}